\documentclass[11pt,english]{article}
\usepackage[T1]{fontenc}
\usepackage[latin9]{inputenc}
\usepackage{geometry}
\usepackage{graphicx}
\usepackage{verbatim}
\usepackage{amsmath}
\usepackage{amssymb}

\makeatletter
\usepackage{times}
\usepackage{algorithm}
\usepackage{algorithmic}
\usepackage{amsthm}
\usepackage{psfrag}
\usepackage[all]{xy}

\makeatletter
 \theoremstyle{plain}
\newtheorem{thm}{Theorem}[section]
  \theoremstyle{plain}
  
 \theoremstyle{definition}
  
  \theoremstyle{definition}
  
  \theoremstyle{definition}
  
  \theoremstyle{definition}
  
  \theoremstyle{plain}
  
  \theoremstyle{plain}
  
    \theoremstyle{plain}
  
     \theoremstyle{plain}
  \newtheorem{conj}[thm]{Conjecture}
     \theoremstyle{plain}     
  
     \theoremstyle{definition}

 \makeatother

\renewcommand{\epsilon}{\varepsilon}

\newcommand{\bu}{\textbf{u}}
\newcommand{\be}{\textbf{e}}

\usepackage{babel}
\makeatother

\title{The Fiedler Rose: On the extreme points of the Fiedler vector}

\author{Lawrence Christopher Evans\footnote{Department of Mathematics, University of Missouri-Columbia. E-mail: evanslc@missouri.edu}}

\begin{document}

\maketitle

\abstract{In this paper I present a counter-example to the conjecture: The Fiedler vector for the graph Laplacian of a tree has its most extreme values at the verticies which are the farthest apart. This counter-example looks roughly like a flower and so I have named it the ``Fiedler rose''.}

\section{Introduction}

Given a simple graph $G=(V,E)$ with vertex set $V=\{v_i:\ i=1,\ldots, n\}$ and edge set $E=\{e_{ij}:\ i\neq j,\ i,j=1,\ldots, n\}$, its \emph{adjacency matrix} is the matrix
$$
A_{ij} = 
\begin{cases}
1,\ \ i\neq j, e_{ij}\in E\\
0,\ \ \text{otherwise},
\end{cases}
$$
its \emph{degree matrix} is the diagonal matrix $D$, where 
$$
D_{ii}=\text{deg}(v_i)
$$
and its \emph{graph Laplacian} is the matrix $L=A-D$.

The graph Laplacian is the discrete analog of the usual Laplacian; for example the heat equation $\frac{du}{dt}=\Delta u$ can be translated to the discrete setting as the system of ODE 
\begin{equation}
\label{DiscreteHeat}
\frac{d\bu}{dt}=L\bu, \bu(0)=\bu_0,
\end{equation}
where the vector $\bu(t)$ denotes the heat at each vertex at time $t$.

Many of the standard results about the Laplacian carry over to the graph Laplacian (and indeed are easier to prove in the discrete setting!): $-L$ is a symmetric, positive-definite matrix, and as such has real non-negative eigenvalues $0=\lambda_1<\lambda_2\leq\lambda_3\leq \ldots \leq \lambda_n$. $\lambda_1=0$ corresponds to the eigenvector whose entries are all the same; that $0$ is an eigenvalue of multiplicity one is a consequence of the Perron-Frobenius theorem. 

Therefore the first eigenvalue of interest is $\lambda_2$. This eigenvalue was first studied extensively by Miroslav Fiedler and he referred to it as the \emph{algebraic connectivity} of the graph due to its connection to connectivity properties of the graph (see~\cite{FiedlerA} and~\cite{FiedlerE}). In honor of Fiedler, it's associated eigenvector has come to be known as the \emph{Fiedler vector} (In general the algebraic connectivity may correspond to an eigenspace of dimension greater than 1, in which case there are many Fiedler vectors. In the example I will consider, however, this is not an issue).

The Fiedler vector has an important role in spectral graph theory due to its successful application towards the problem of graph partitioning (see~\cite{Partition}). But what motivates the example I will present is its relation to the discrete heat equation.

The solution to (\ref{DiscreteHeat}) is given by
$$
\bu(t)=\sum_{i=1}^n (\bu_0,\be_i)e^{-\lambda_i t}\be_i,
$$
where the $\be_i$ are the orthonormal basis of eigenvectors of $L$. As $\lambda_i>0$ for $i\geq 2$ and $\lambda_1=0$ corresponds to $\be_1=\frac{1}{\sqrt{n}}(1,1,\ldots,1)^T$, we have that $\bu(t)\to(\bu_0,\be_1)\be_1$. That is, the heat will eventually even out until, in the limit, the heat is constant across all verticies.

To study the long term behavior of $\bu(t)$, we note that for $t$ large, assuming $\lambda_2\neq\lambda_3$ and $(\bu_0,\be_2)\neq 0$, 
$$
\bu(t)\approx(\bu_0,\be_1)\be_1+(\bu_0,\be_2)e^{-\lambda_2 t}\be_2,
$$
as the other terms die out faster. That is, in the long run, $\bu(t)$ has the same structure as the Fiedler vector $\be_2$ up to a constant multiple and translation. Therefore, the Fiedler vector captures the transient behavior of the heat flow. In particular, the extreme points of $\bu(t)$ will be extreme points of the Fiedler vector.

Therefore, heuristically, the extreme points of the Fiedler vector should correspond to the most ``insulated verticies'', since these are the verticies in which heat/cold will get trapped and $\bu(t)$ will stay the hottest/coldest. On a tree, it might seem at first that the two most insulated verticies would be the ones farthest apart. This raises the following conjecture.
\begin{conj}
\label{TreeConjecture}
Suppose $G=(V,E)$ is a tree. If $v^*,w^*\in V$ are such that
$$
|e_2(v)-e_2(w)|\leq |e_2(v^*)-e_2(w^*)|\quad \forall\ (v,w)\in V^2,
$$
then
$$
d(v,w)\leq d(v^*,w^*)\quad \forall\ (v,w)\in V^2,
$$
where $d(v,w)$ is the graph-distance between $v$ and $w$. In other words, the extreme values of the Fiedler vector are among the pairs of verticies which are the farthest apart.
\end{conj}

Indeed, in a recent paper by Chung, Seo, Adluru, and Vorperian,~\cite{Conjecture}, a similar conjecture is made (see Section \ref{LastSection} for the precise statement of their conjecture).

However, I will now present a counter-example to Conjecture \ref{TreeConjecture}. Consider the graph in Figure 2.

\begin{figure}[!htbp] \label{Figure1}\caption{The Fiedler rose}\centering \includegraphics[scale=0.4]{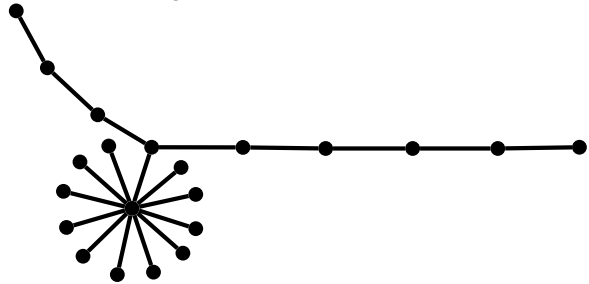}
\end{figure}

If we consider the Fiedler vector (which will be computed in the next section) for this graph and color the verticies red, whose entry in the Fiedler vector is positive, and blue, whose entry in the Fiedler vector is negative, we get Figure 2.

\begin{figure}[!htbp] \label{Figure2}\caption{Heat map of the Fiedler vector}\centering \includegraphics[scale=0.4]{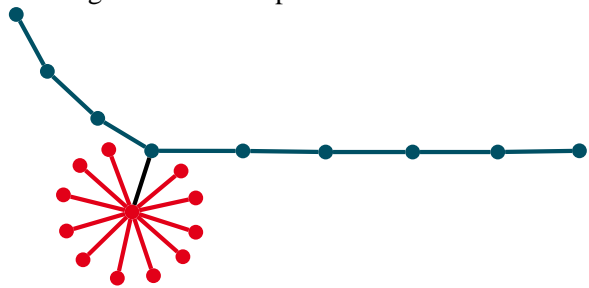}
\end{figure}

The graph now looks like a rose with a curved ``leaf'', long ``stem'', and many ``petals'', hence I will refer to it as the \emph{Fiedler rose}. Note that the extreme values of the Fiedler vector are not the verticies farthest apart, the tip of the leaf and stem, but are instead the stem and any of the petal verticies (the exact values of the Fiedler vector will be computed in the next section). Thus, the Fiedler rose shows false Conjecture \ref{TreeConjecture}.

\section{Computation}

I first label the verticies of the rose as depicted, letting $s$ be the number of verticies in the stem and $p$ be the number of petal-verticies (I will consider Fiedler roses of varying stem lengths and petal counts). 

\begin{figure}[!htbp] 
\centering 
\label{Figure3}
\includegraphics[scale=0.4]{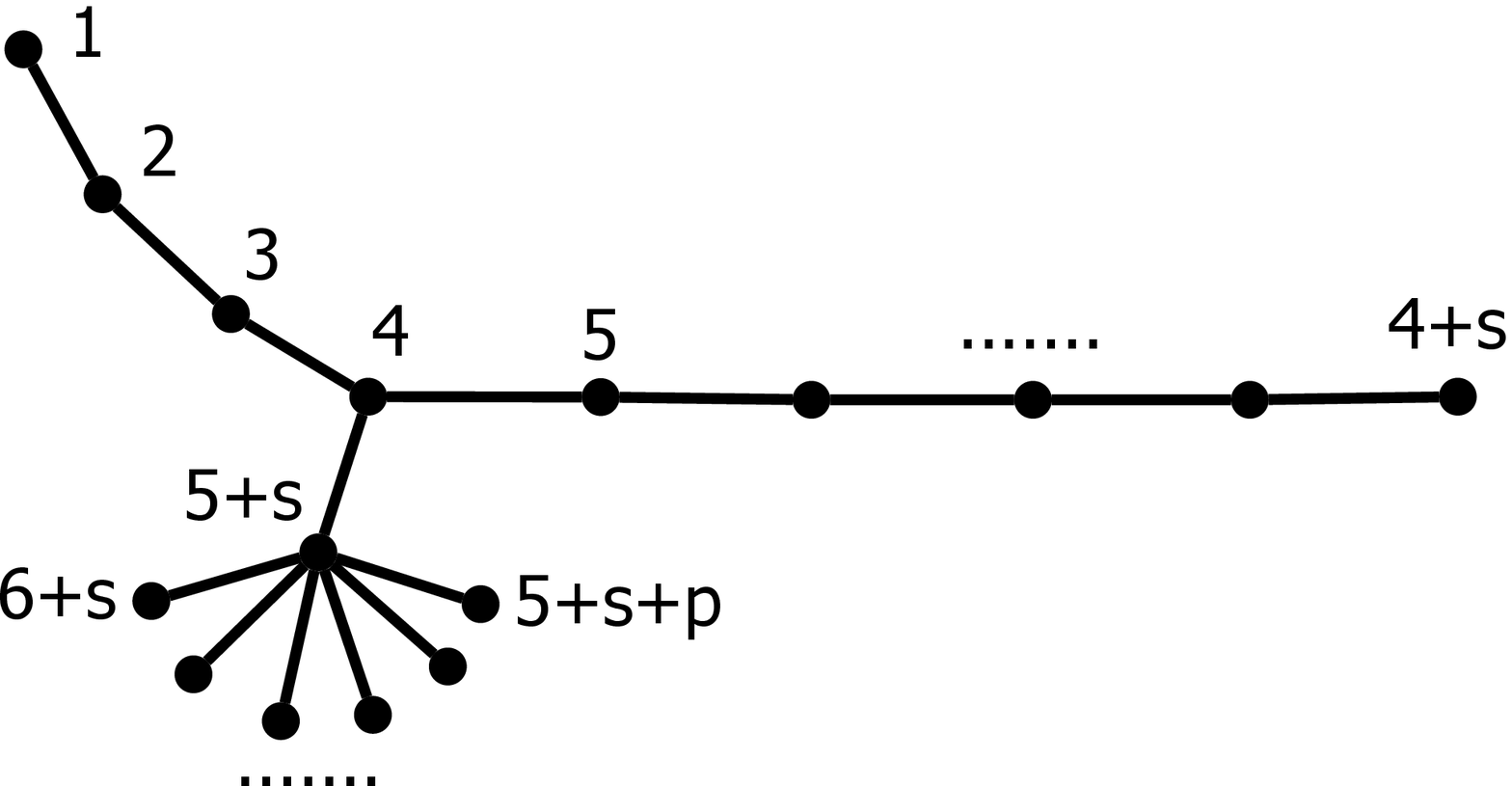}
\end{figure}

I then use the following MATLAB code to compute the Fiedler vector, which displays the Fiedler vector by embedding it in the matrix $B$.

\scriptsize 
\begin{verbatim}
>function FiedlerRose(p,s)
>
>N=5+p+s; %Total verticies = 5 for leaf/rose center + s for stem + n petals
>L=zeros(N,N);
>r=5+s; %Center rose node
>
>for i=1:(4+s-1),  %Building the leaf and stem
>    j=i+1;
>    L(i,i)=L(i,i)-1;
>    L(i,j)=L(i,j)+1;
>    L(j,i)=L(j,i)+1;
>    L(j,j)=L(j,j)-1;
>end
>
>i=4; %Connecting the stalk to the rose center
>j=r;
>L(i,i)=L(i,i)-1;
>L(i,j)=L(i,j)+1;
>L(j,i)=L(j,i)+1;
>L(j,j)=L(j,j)-1;
>
>for j=r+1:r+p, %Building the p rose petals
>    i=r;
>    L(i,i)=L(i,i)-1;
>    L(i,j)=L(i,j)+1;
>    L(j,i)=L(j,i)+1;
>    L(j,j)=L(j,j)-1;
>end
>
>[V,D]=eigs(L,3,0.001); %Grabs the 3 eigenvalues/vectors with eigenvalue closest
>                       %to 0.001.
>v=V(:,2); %We set v to be the Fiedler vector.
>if v==V(:,3), end %Abort if there is more than one Fiedler vector.
>if v(r)<0; %We multiply by -1 if neccesary to ensure the center of rose is 'hot'.
>   v=-1*v; 
>end
>
>B=zeros(2,4+s); %We initialize B where we will store the rose for display.
>index=0; %The index keeps track of which vertex we are currently placing.
>
>for j=1:4+s;
>    index = index+1;
>    B(1,j)=v(index); %Place vertex into appropriate spot in B.
>end
>index = index+1;
>B(2,4) = v(index); %Place center of rose.
>index = index+1;
>B(2,3) = v(index); %Place petal of rose. 
>                   %All other petals have the same value by symmetry
>                   %so it suffices to display one.
>B
\end{verbatim}
\normalsize

If we run the code with  $p=11$ and $s=5$ (the values for the graph we considered in the previous section), we get the following output:

\scriptsize 
\begin{verbatim}
B =

   -0.0093   -0.0085   -0.0071   -0.0051   -0.1481   -0.2793   -0.3881   -0.4659   -0.5064
         0         0    0.1525    0.1403         0         0         0         0         0
         
\end{verbatim}
\normalsize
Here, the top row of $B$ consists of the values of the Fiedler vector along the leaf and stem, where the left-most value is the value at the tip of the leaf and the right-most value is the value at the tip of the stem. The entry $B(2,4)$ is the value of the Fiedler vector at the center of the rose, and $B(2,3)$ is the value at any of the petals of the rose, which are all the same by symmetry (The matrix $B$ preserves the overall shape of the rose for easy viewing).

As was claimed earlier, the largest value of the Fiedler vector is at a petal vertex and the smallest value is at the tip of the stem.

If we remove a petal and run the code with $p=10$ and $s=5$, however, we get a much different output:

\scriptsize 
\begin{verbatim}
B =

    0.0074    0.0068    0.0056    0.0040   -0.1414   -0.2752   -0.3865   -0.4662   -0.5077
         0         0    0.1606    0.1474         0         0         0         0         0
         
\end{verbatim}
\normalsize

Now the leaf of the rose is also positively valued (and if we were to color the verticies by their sign the picture would no longer resemble a rose -- the leaf would be red as well!). 

As we decrease the number of petals further, the value of the Fiedler vector at the tip of the leaf rises until it eventually becomes the most positive value. Indeed, if we run the code with $p=3$ and $s=5$, we see that the most extremely valued verticies are the tip of the leaf and stem:

\scriptsize 
\begin{verbatim}
B =

    0.2514    0.2253    0.1758    0.1081   -0.0597   -0.2213   -0.3600   -0.4612   -0.5147
         0         0    0.2198    0.1970         0         0         0         0         0
         
\end{verbatim}
\normalsize

\section{Concluding Remarks}
\label{LastSection}

First, it should be noted that the actual conjecture posed as Conjecture 2 in~\cite{Conjecture} is that the second eigenfunction of the Laplacian on a smooth manifold without boundary achieves its extreme values at the points which are the furthest geodesic distance apart. The authors then go on to state that ``Conjecture 2
can be applicable not only to differentiable manifolds but to graphs and surface
meshes as well'' The Fiedler rose shows that this conjecture is false for graphs, even if one assumes the graph is a tree.

It may be possible to turn the Fiedler rose into a counter-example for Conjecture 2 of~\cite{Conjecture} by ``puffing up'' each of the edges and considering the surface of the resulting 3-dimensional object... but I do not know an easy way to computationally verify this.


\bibliography{Bibliography}{}

\newcommand{\noopsort}[1]{} \newcommand{\printfirst}[2]{#1}
  \newcommand{\singleletter}[1]{#1} \newcommand{\switchargs}[2]{#2#1}
  \def\cprime{$'$}
\begin{thebibliography}{1}

\bibitem{Conjecture}
Moo Chung, Seongho Seo, Nagesh Adluru, and Houri Vorperian.
\newblock Hot spots conjecture and its application to modeling tubular
  structures.
\newblock In Kenji Suzuki, Fei Wang, Dinggang Shen, and Pingkun Yan, editors,
  {\em Machine Learning in Medical Imaging}, volume 7009 of {\em Lecture Notes
  in Computer Science}, pages 225--232. Springer Berlin / Heidelberg, 2011.

\bibitem{FiedlerA}
Miroslav Fiedler.
\newblock Algebraic connectivity of graphs.
\newblock {\em Czechoslovak Math. J.}, 23(98):298--305, 1973.

\bibitem{FiedlerE}
Miroslav Fiedler.
\newblock A property of eigenvectors of nonnegative symmetric matrices and its
  application to graph theory.
\newblock {\em Czechoslovak Math. J.}, 25(100)(4):619--633, 1975.

\bibitem{Partition}
Alex Pothen, Horst~D. Simon, and Kang-Pu Liou.
\newblock Partitioning sparse matrices with eigenvectors of graphs.
\newblock {\em SIAM J. Matrix Anal. Appl.}, 11(3):430--452, 1990.
\newblock Sparse matrices (Gleneden Beach, OR, 1989).

\end{thebibliography}
\bibliographystyle{plain}

\end{document}